\documentclass{IEEEtran4PSCC}

\usepackage{xcolor}
\usepackage[normalem]{ulem}
\ifCLASSINFOpdf
   \usepackage[pdftex]{graphicx}
\else
   \usepackage[dvips]{graphicx}
\fi
\usepackage{comment}
\usepackage{optidef}
\usepackage[]{amsmath}

\usepackage{tablefootnote}

\usepackage{algorithm}
\usepackage{algpseudocode}
\usepackage{dblfloatfix}
\usepackage{multirow}
\usepackage[caption=false,font=footnotesize]{subfig}
\makeatletter
\let\old@ps@headings\ps@headings
\let\old@ps@IEEEtitlepagestyle\ps@IEEEtitlepagestyle
\def\psccfooter#1{%
    \def\ps@headings{%
        \old@ps@headings%
        \def\@oddfoot{\strut\hfill#1\hfill\strut}%
        \def\@evenfoot{\strut\hfill#1\hfill\strut}%
    }%
    \def\ps@IEEEtitlepagestyle{%
        \old@ps@IEEEtitlepagestyle%
        \def\@oddfoot{\strut\hfill#1\hfill\strut}%
        \def\@evenfoot{\strut\hfill#1\hfill\strut}%
    }%
    \ps@headings%
}
\makeatother

\psccfooter{%
        \parbox{\textwidth}{\hrulefill \\ \small{22nd Power Systems Computation Conference} \hfill \begin{minipage}{0.2\textwidth}\centering \vspace*{4pt} \includegraphics[scale=0.06]{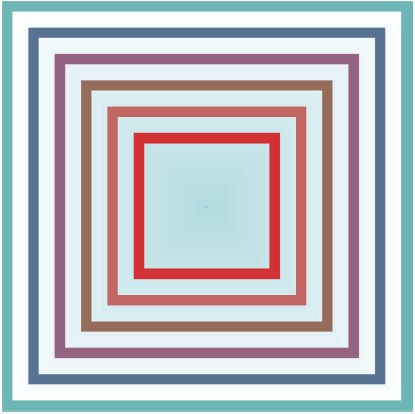}\\\small{PSCC 2022} \end{minipage} \hfill \small{Porto, Portugal --- June 27 -- July 1, 2022}}%
}

\begin{document}
\newcommand{\squeezeuptwo}{\vspace{-2mm}}
\newcommand{\squeezeupthree}{\vspace{-3mm}}
\title{Three-Phase Infeasibility Analysis for Distribution Grid Studies}

\author{
\IEEEauthorblockN{Elizabeth Foster, Amritanshu Pandey, and Larry Pileggi}
\IEEEauthorblockA{ Dept. of Electrical and Computer Engineering \\
Carnegie Mellon University\\
Pittsburgh, PA USA\\
\{emfoster, amritanp, pileggi\}@andrew.cmu.edu}
}


\maketitle


\begin{abstract}
With the increase of distributed energy resources in the distribution grid, planning to ensure sufficient infrastructure and resources becomes critical. Planning at the distribution level is limited by the complexities of optimizing unbalanced systems. In this paper we develop a three-phase infeasibility analysis that identifies weak locations in a distribution network. This optimization is formulated by adding slack current sources at nodes in the system and minimizing their norm subject to distribution power flow constraints. Through this analysis we solve instances of power flow that would otherwise be infeasible and diverge. Under conditions when power flow is feasible, our approach is equivalent to standard three-phase power flow; however, for cases where power flow fails, the nonzero slack injection currents compensate for missing power to make the grid feasible. Since an uncountable number of injected currents can provide feasibility, we further explore the optimization formulation that best fits the solution objective through use of both a least squares and an L1 norm objective. Our L1 norm formulation localizes power deficient locations through its inherent sparsity. We show the efficacy of this approach on realistic unbalanced testcases up to 8500 nodes and for a scenario with a high penetration of electric vehicles.
\end{abstract}


\thanksto{\noindent Submitted to the 22nd Power Systems Computation Conference (PSCC 2022).}

\section{Introduction}
An increasingly complex distribution grid with more distributed energy resources (DERs) and inverter-based resources necessitates more proactive planning to ensure there is adequate infrastructure and resources for changing demand \cite{NAS}. Planning at the distribution level is often reactive because there is a limited offering of tools capable of performing three-phase steady-state optimizations and utilities themselves may lack high fidelity data at the distribution level \cite{PNNLDER}. As a greater number of sensors and smart meters are introduced on the distribution grid and as utilities collect more detailed information on the distribution grid layout, one of the impediments to more robust and proactive distribution analysis is the lack of three-phase optimization methods for large-scale networks. 
\begin{figure*}[btp]
    \centering
    \includegraphics[width=\textwidth]{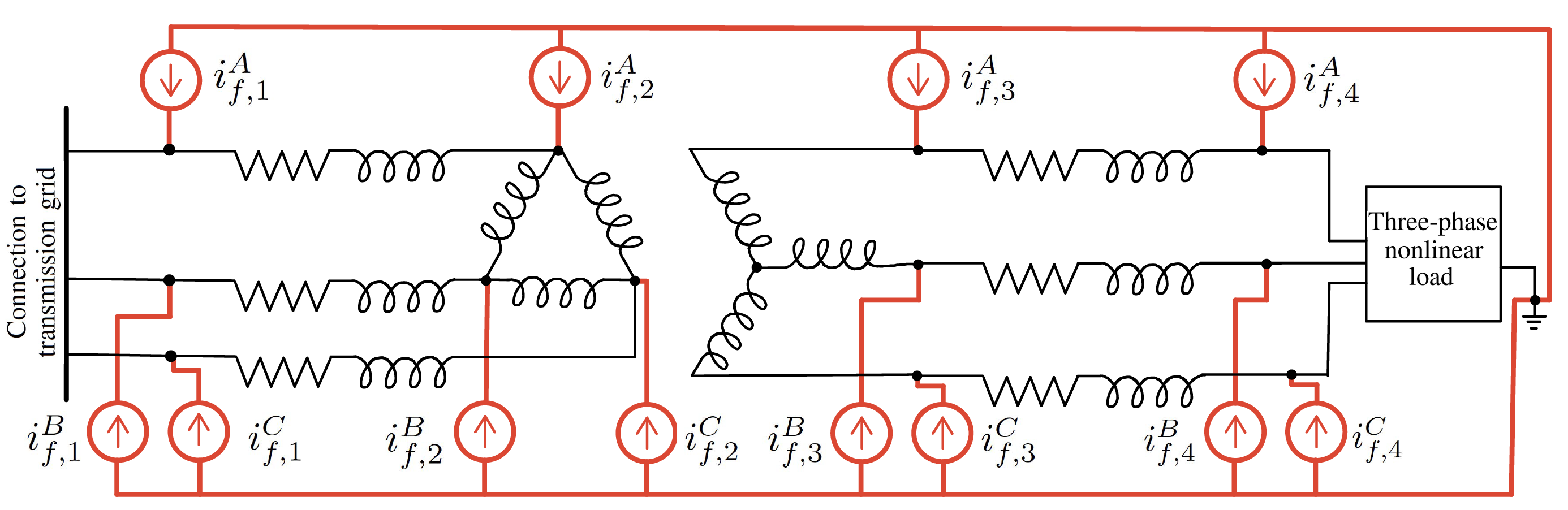}
    \caption{Representation of the IEEE 4 node case with infeasibility current sources at each node}
    \label{4buscase}
\end{figure*}
When considering planning for future load growth at the grid-edge, simulating three-phase power flow is not enough to provide advanced planning insight. Power flow simulations can fail because the testcase misrepresents a real physical system or because no feasible solution exists for the given set of operating conditions. It may also fail due to choice of poor initial conditions or unavailability of good initial conditions \cite{Amrit}. System planners require clear indications regarding why power flow simulations fail since different causes correspond to different remedies as it relates to planning capabilities. We lack methods that can not only identify an infeasible or collapsed network, but also provide useful insight on the failure. 

There is some prior work in identifying infeasible grid settings in both positive sequence networks \cite{Overbye} \cite{Marko} \cite{Cindy} and distribution networks \cite{Dall}. \cite{Overbye} introduced the concept of \textit{missing power} at the transmission level based on power mismatch equations augmented with slack variables. \cite{Dall} introduced a similar concept for distribution networks to identify problematic optimal power flow constraints. More recently, we developed a more robust approach using equivalent circuit formulation (ECF) for transmission systems in \cite{Marko}. Our work in \cite{Cindy} built on \cite{Marko} by adding an L1 norm regularization to induce greater sparsity in the solution vector. Nonetheless, none of these methods are robust, scalable, and localized for three-phase distribution networks.
Unbalance between phases is common in distribution networks, so positive sequence methods in \cite{Overbye} \cite{Marko} \cite{Cindy} are not directly applicable because they lack the necessary but complicated phase relationships. \cite{Overbye}\cite{Dall} do not demonstrate scalability beyond a few hundred nodes. \cite{Overbye} \cite{Marko} \cite{Dall} do not attempt to localize the sources of infeasibility. Difficulties in three-phase distribution systems further arise from a low X-R ratio; significant number of direct, non-controllable customers; and significant levels of unbalance \cite{Gerth}. It is an imminent need in distribution planning studies to have large scale three-phase analysis capable of detecting localized sources of infeasibility, and none of the existing works in this space can fill the void. 

Therefore, to address this research gap, we propose an ECF-based three-phase optimization method called the three-phase infeasibility analysis (TPIA) that can quantify and localize weak locations in the distribution grid. TPIA can inform system planners where new or updated infrastructure may be needed and it helps grid operators locate areas in the system where reliability issues may occur under abnormal or extreme conditions. Section II details relevant prior work. 

In TPIA, we first construct aggregated equivalent circuits to represent three-phase distribution grids. To solve otherwise infeasible networks, we add variable slack current sources at each node in the constructed circuit. These sources then compensate for the mismatch in the Kirchhoff’s Current Law (KCL) equations. Due to the discrepancy between the number of unknown variables and constraints, the solution of this overall circuit will have infinitely many solutions. Because of this mismatch, we minimize the norm of the slack sources in TPIA through an optimization approach. 

An illustration of the slack current source injection is shown for the IEEE 4 node testcase in Fig. \ref{4buscase}. The slack current sources, which we call ``infeasibility currents" are the red independent current sources labeled $i_f$. Within the optimization, these new infeasibility currents are minimized and constrained by KCL at each node. In this paper, we formulate TPIA with two different objective functions. We first use the common least squares objective and then we develop TPIA using the least absolute value objective (L1 norm). We use the L1 norm to induce greater sparsity in results in comparison to the least squares approach and to localize sources of a power flow simulation failure. Methodology for the two formulations is shown in Section III. Implementation of TPIA is shown in Section IV. 

Section V describes experimental set-ups and results for three-phase networks up to 8500 nodes and discusses the empirical differences between the two TPIA formulations. We show that both the least squares and L1 norm formulations of TPIA converge when three-phase power flow simulations in GridLAB-D and SUGAR-D (from our work in \cite{Amrit}) otherwise fail. We show that the L1 formulation localizes grid infeasibilities whereas the least squares formulation does not. Finally, we provide a use case with the L1 TPIA formulation for evaluating electric vehicle (EVs) penetration and validate our results by constructing a feasible network based on information from the L1 TPIA results.

\section{Prior Work}

TPIA solves an optimization problem to locate weak-spots in a three-phase unbalanced distribution network especially when power flow solution does not exist. No standard tool today (e.g. OpenDSS, GridLAB-D) can identify these weak spots for infeasible large-scale realistic three-phase networks. We develop the TPIA algorithm building on state-of-the-art research in three-phase distribution grid analysis. We briefly discuss the work here.


\subsection{Steady-state Analysis of Distribution Grid}

Distribution grid analysis can be fundamentally categorized as either a simulation (or power flow) or an optimization problem.

\subsubsection{Three-Phase Power Flow}
Three-phase power flow has been formulated using a number of techniques like the Backward-Forward Sweep Method (BFS) \cite{BFS}, Current Injection Method (CIM) \cite{PAN}, and Equivalent Circuit Formulation (ECF) \cite{Amrit}. BFS works well for weakly-meshed or radial systems but is known to suffer otherwise \cite{BFS}. CIM improves solution robustness by formulating the problem as current-balance rather than power-balance \cite{PAN}. However, CIM is known to suffer when there is a high penetration of generation models (PV nodes), which is increasingly common in distribution grids. ECF uses circuit theory to model the electric grid as a three-phase circuit and uses circuit-theoretic heuristics to obtain the solution to the three-phase power flow problem. We build the TPIA framework with ECF as the basis. 

In ECF, the AC network constraints are written using KCL, which requires that the sum of currents at each node add up to zero. Mathematically, the formulation is shown in \eqref{eqn: PowerFlow-1}-\eqref{eqn: PowerFlow-4}, where $\odot$ and $\oslash$ represent element-wise multiplication (Hadamard Product) and element-wise division (Hadamard Division) respectively. Variables
$V_R$ and $V_I$ are vectors of $n$ unknown variables in rectangular components; $\Omega$ represents the phases $(\Omega={A,B,C})$; $G$ and $B$ are respectively the admittance and susceptance matrices of size $n \times n$; and $I^{NL}_R$ and $I^{NL}_I$ are $n \times 1 $ vectors representing the nonlinear real and imaginary load or generator currents. $I^{NL}_R$ and $I^{NL}_I$ are functions of unknown variables and are shown in \eqref{eqn: PowerFlow-3} and \eqref{eqn: PowerFlow-4} where $P^{\Omega}$ represents real power and $Q^{\Omega}$ represents reactive power. For a PV bus, $Q^{\Omega}$ is an unknown variable. At the PQ bus, $Q^{\Omega}$ is known so the current equations can be simplified. While the ECF approach can handle PV buses, for brevity, \eqref{eqn: PowerFlow-1}-\eqref{eqn: PowerFlow-4} do not include additional PV model terms and constraints. 
\begin{align}
    G V^{\Omega}_R - B V^{\Omega}_I + I^{NL}_R(V^\Omega_R, &V^\Omega_I, Q^{\Omega}) = 0 \label{eqn: PowerFlow-1} \\
    B V^{\Omega}_R + G V^{\Omega}_I + I^{NL}_I(V^\Omega_R, &V^\Omega_I, Q^{\Omega}) = 0 
    \label{eqn: PowerFlow-2} \\
   I^{NL}_R(V^{\Omega}_R, V^{\Omega}_I, Q^{\Omega}) = P^{\Omega}\odot &V_R^{\Omega}  + Q^{\Omega}\odot V_I^{\Omega} \nonumber\\ 
   & \oslash ((V_R^{\Omega})^2 + (V_I^{\Omega})^2)  \label{eqn: PowerFlow-3} \\
    I^{NL}_I(V^{\Omega}_R, V^{\Omega}_I, Q^{\Omega}) = P^{\Omega}\odot &V_I^{\Omega}  - Q^{\Omega}\odot V_R^{\Omega} \nonumber \\
   & \oslash ((V_R^{\Omega})^2 + (V_I^{\Omega})^2)   \label{eqn: PowerFlow-4}
\end{align}
 


\subsubsection{Optimizations in Three-Phase Distribution Grids}

While many earlier works on distribution grid optimization relied on linearized DC constraints that poorly represent the system \cite{Gerth}, there has been significant advances in AC-network constrained optimization methods for distribution grids. For instance,  \cite{Dall} maps the AC network constraints using power balance equations and \cite{Oconnell} uses CIM to map the AC network constraints with KCL. Both of these approaches can work with general objectives. However, because these existing distribution grid optimization formulations do not demonstrate robustness and scalability to solve large-scale optimization problems ($>$5000 nodes), we extend the ECF framework to include AC-network constrained optimization for use with TPIA.

\subsection{Optimization Techniques for Identifying Infeasible Power Grids}

There are some approaches that can identify infeasible positive-sequence transmission networks \cite{Overbye}, \cite{Marko}, \cite{Cindy} as well as distribution networks \cite{Dall}. For transmission positive-sequence networks, \cite{Overbye}, \cite{Marko}, \cite{Cindy} have all developed techniques to detect grid infeasibilities. \cite{Overbye}'s approach is based on power-mismatch equations. This approach is known to struggle with three-phase distribution networks due to the use of power balance equations \cite{Amrit}.  Another distinct approach for identifying grid infeasibility was introduced for transmission positive sequence networks using a circuit-theoretic ECF foundation \cite{Marko}. However, both \cite{Overbye} and \cite{Marko} cannot localize the sources of infeasibility in the grid. Our follow-on work in \cite{Cindy} localized grid infeasibilities in transmission networks by formulating the L2 norm objective function with an L1 norm regularization term. This formulation encourages some sparsity in the optimization solution, which can be used to find regions of grid infeasibility. However, none of these works apply directly to three-phase distribution networks due to its unbalanced and complex nature. 

Research in identifying infeasible three-phase distribution networks is sparse but exists. For instance, \cite{Dall} formulates a three-phase optimization problem using power balance equations to identify infeasible distribution grids. It adds and minimizes slack variables to accommodate constraint violations and identify infeasible networks. However, many challenges persist with this approach. The slack variables do not have a physical meaning, so their values provide no grid insight outside of a constraint violation; the method has not been shown to scale; and the approach is unable to localize sources of infeasibilities. With TPIA, we will address all these challenges and develop an infeasibility analysis framework that provides practical, interpretable causes of distribution grid infeasibility.

\section{Three-Phase Infeasibility Analysis}
The goal of TPIA is to identify and localize weak sections of three-phase networks when the three-phase power flow (TPF) problem has no solution. Thus we begin by constructing a three-phase equivalent circuit of the evaluated distribution network (compare to \cite{Amrit}) and transforming the TPF simulation problem into a TPIA optimization problem.

To illustrate further, consider the 4-bus network in Fig. \ref{4buscase}. Without the presence of current sources in red, the network can be solved directly using the ECF formulation in \cite{Amrit} to obtain a TPF solution. However, if there is no physical solution for the problem, the traditional TPF technique would diverge even with the best heuristics. Instead, in TPIA, we add variable current sources at each phase of each node, namely the infeasibility sources $i_{f,k}^{\Omega}$ (shown in red for node $k$), to ensure solvability of the network. The goal of the variable infeasibility sources is to supply any missing current where the KCL equations cannot be satisfied within the confines of the original network. While the infeasibility sources ensure feasibility of the circuit (and correspondingly the TPF problem) for all possible operating conditions, infinitely many solutions exist as the system is underdetermined with these additional variable $i_{f}^{\Omega}$ sources. To ensure that TPIA returns either the TPF solution (when the problem is feasible) or the least number of infeasibility currents (when the problem is infeasible), we minimize the norm of $i_{f}^{\Omega}$ subject to TPF AC network constraints.

\subsection{General Formulation}
 Mathematically, TPIA takes the form in \eqref{eqn: generic}, where $i_{f}^{\Omega}$ represents a vector of infeasibility currents; $\mathcal{G}(i_{f}^{\Omega})$ is a generic objective function; and $\mathcal{H}(X)$ represents three phase AC network constraints from CIM or ECF where $X$ is the vector of unknown variables in TPF. Unlike positive-sequence power flow, $\mathcal{H}(X)$ models each phase of $\Omega$ in the distribution network and can capture both balanced and unbalanced networks. 

\begin{mini}|s|
{X, i_{f}^{\Omega}}{\mathcal{G}(i_{f}^{\Omega})}
{}{}
\addConstraint{\mathcal{H}(X) - i_{f}^{\Omega} = 0}
\label{eqn: generic}
\end{mini}

In this paper, we derive and solve two distinct forms of the objective functions $\mathcal{G}(i_{f}^{\Omega})$ in \eqref{eqn: generic}. The first minimizes the  infeasibility currents using least squares and the second minimizes using the L1 norm.

\subsection{Location of Infeasibility Sources}
Infeasibility current sources can be added at any defined set of nodes. The ability to optimize over a predefined set of nodes is beneficial because it allows a system planner to identify which locations in the system they are capable of modifying before conducting the optimization. The L1 norm approach inherently returns sparse solutions which effectively localize weak spots in the distribution grid whereas the least squares approach spreads out infeasibility currents. Therefore, defining a subset of notes is especially beneficial in the least squares approach. Depending on available resources, system planners may prefer either to make larger adjustments at a few nodes based on the localized results from L1 norm approach or to make a greater number of smaller adjustments via the least squares solution. While the inherent sparsity of the L1 norm is an advantage, it is non-differentiable. The standard least squares approach is empirically much easier and faster to solve in comparison. We explore both formulations in the following subsections because each has distinct advantages.

\subsection{Least Squares Formulation}
As a first approach to TPIA, we use least squares to represent the objective function in \eqref{eqn: generic}. The least squares objective function is well studied within optimization and adds no inequality constraints to the Lagrangian formulation as it has a continuous derivative. Reformulating \eqref{eqn: generic} with a least-squares objective function and using the TPF equations in \eqref{eqn: PowerFlow-1}-\eqref{eqn: PowerFlow-4}, we obtain the formulation in \eqref{eqn: opt1} - \eqref{eqn: opt3}.
\begin{subequations}
\begin{IEEEeqnarray}{s,rCl'rCl'rCl}
min   & \IEEEeqnarraymulticol{9}{l}{\frac{1}{2}||i_{f,R}^{\Omega}||^2+\frac{1}{2}||i_{f,I}^{\Omega}||^2} \label{eqn: opt1} \\
s.t. & G V^{\Omega}_R - B V^{\Omega}_I + I^{NL}_{R}(V^{\Omega}_R, V^{\Omega}_I, Q^{\Omega}) - i_{f,R}^{\Omega} = 0 
\label{eqn: opt2}\\
           & B V^{\Omega}_R + G V^{\Omega}_I + I^{NL}_{I}(V^{\Omega}_R, V^{\Omega}_I, Q^{\Omega}) - i_{f,I}^{\Omega} = 0 \label{eqn: opt3}
\end{IEEEeqnarray}
\end{subequations}
To find an optimal solution to this problem, we formulate the Lagrangian of the problem in \eqref{eqn: L2lagrangian} where $\lambda_R^{\Omega}$ and $\lambda_I^{\Omega}$ are the vector of dual variables or KKT multipliers. In least squares based TPIA, we find that these dual variables are closely tied to the infeasibility currents, such that non-zero magnitudes for dual variables $(|\lambda_k^\Omega|>0)$ at node $k$ exist only when the corresponding node $k$ has non-zero infeasibility currents, $|i_{f,k}^{\Omega}|>0$.

\begin{align}
    \mathcal{L}(i_{f,R}^{\Omega}, &i_{f,I}^{\Omega}, V_R^{\Omega},  V_I^{\Omega}, Q^{\Omega}, \lambda_R^{\Omega}, \lambda_I^{\Omega}) =  \label{eqn: L2lagrangian}\\ 
    & \frac{1}{2}||i_{f,R}^{\Omega}||^2+\frac{1}{2}||i_{f,I}^{\Omega}||^2 \nonumber\\
    + & \lambda_R^{\Omega,T}\left(GV^{\Omega}_R - BV^{\Omega}_I + I^{NL}_{R}(V^{\Omega}_R, V^{\Omega}_I, Q^{\Omega}) - i_{f,R}^{\Omega} \right) \nonumber\\
    + & \lambda_I^{\Omega,T} \left(B V^{\Omega}_R + G V^{\Omega}_I + I^{NL}_{I}(V^{\Omega}_R, V^{\Omega}_I, Q^{\Omega}) - i_{f,I}^{\Omega} \right) \nonumber
\end{align}

With the Lagrangian formulated, we derive the first-order optimality conditions for the KKT conditions. The KKT conditions for the primal variables are shown in Equations \eqref{L2statIR} - \eqref{L2statVI} and are solved in our optimization routine using Newton's Method.

\begin{subequations}
\begin{IEEEeqnarray}{ll}
    \frac{\partial \mathcal{L}}{\partial i_{f,R}^{\Omega}}  = & i_{f,R}^{\Omega} - \lambda_R^{\Omega} = 0 \label{L2statIR} \\
    \frac{\partial \mathcal{L}}{\partial i_{f,I}^{\Omega}}  = & i_{f,I}^{\Omega} - \lambda_I^{\Omega} = 0  \label{L2statII}\\
    \frac{\partial \mathcal{L}}{\partial Q^{\Omega}}  = &  \frac{\partial I^{NL}_R}{\partial Q^{\Omega}}\lambda_R^{\Omega} + \frac{\partial I^{NL}_I}{\partial Q^{\Omega}}\lambda_I^{\Omega} = 0 \label{L2statQ} \\
    \frac{\partial \mathcal{L}}{\partial V^{\Omega}_R} = &
    (G^T + \frac{\partial I^{NL}_R}{\partial V^{\Omega}_R})\lambda_R^{\Omega} + (B^T + \frac{\partial I^{NL}_I}{\partial V^{\Omega}_R})\lambda_I^{\Omega} = 0 \label{L2statVR}\\
    \frac{\partial \mathcal{L}}{\partial V^{\Omega}_I} = 
     & (-B^T + \frac{\partial I^{NL}_R}{\partial V^{\Omega}_I})\lambda_R^{\Omega} + (G^T + \frac{\partial I^{NL}_I}{\partial V^{\Omega}_I})\lambda_I^{\Omega} = 0 \label{L2statVI}
\end{IEEEeqnarray}
\end{subequations}
In \eqref{L2statIR} and \eqref{L2statII}, we can see that per KKT conditions, $i_{f,R}^{\Omega} = {\lambda}_R^{\Omega} $ and $i_{f,I}^{\Omega} = \lambda_I^{\Omega} $, which validates our earlier observation about the tight coupling between dual variables and infeasibility currents. Equations \eqref{L2statQ} - \eqref{L2statVI} have nonlinear components and will be linearized using a Taylor Series approximation which allows us to use Newton's Method to solve this optimization. The primal feasibility constraints are the TPF equations shown in Equations \eqref{eqn: opt2} - \eqref{eqn: opt3}. 

The disadvantage of the least squares approach is that it does not return a sparse solution vector for infeasibility currents $i_{f}^{\Omega}$. While a system planner would be able to identify an infeasible network with the least squares formulation of TPIA, they would not be able to localize the regions of the network that caused the infeasibility because the infeasibility currents would be spread throughout the system. We show this key characteristic in Fig. \ref{fig_sim}A which demonstrates that infeasibility currents are present at nearly every node in the system during an infeasible operating scenario. Planners are unable to use this information to perform corrective actions. To gain any realistic insight on physical remedies, the planner would need to define a smaller subset of nodes where infeasibility current sources should be added rather than across the whole system. Therefore, to get a more sparse solution vector such that it provides more immediately localized information to the planner, we consider a second approach to TPIA using an L1 norm objective function.


\subsection{L1 Norm Formulation}
We next consider the L1 norm based objective, which can localize the nonzero injection currents to the weakest sections of the grid because of its inherent property to provide sparse solutions. 

With this mechanism, when a network file incorrectly represents the real physical system and TPF fails, system planners can isolate where the modeling error is based on looking at the location of the nonzero injection currents. This will prevent significant manual overhead that would otherwise be required to find the error in the system topology (e.g., typo in the input file). In situations where TPF fails due to inherent lack of resources, the L1 formulation of TPIA can be used as a planning tool to ascertain where new infrastructure may be needed in the grid to accommodate growing or changing electricity demand. 


Mathematically, the L1 norm formulation of TPIA is formulated as: 
\begin{subequations}
\begin{IEEEeqnarray}{s,rCl'rCl'rCl}
min   & \IEEEeqnarraymulticol{9}{l}{|i_{f,R}^{\Omega}|+|i_{f,I}^{\Omega}|} \label{eqn: opt1_L1} \\
s.t. & G V^{\Omega}_R - B V^{\Omega}_I + I^{NL}_{R}(V^{\Omega}_R, V^{\Omega}_I, Q^{\Omega}) - i_{f,R}^{\Omega} = 0 
\label{eqn: opt2_L1}\\
           & B V^{\Omega}_R + G V^{\Omega}_I + I^{NL}_{I}(V^{\Omega}_R, V^{\Omega}_I, Q^{\Omega}) - i_{f,I}^{\Omega} = 0 \label{eqn: opt3_L1}
\end{IEEEeqnarray}
\end{subequations}

Unlike the least squares formulation in \eqref{eqn: opt1} - \eqref{eqn: opt3}, the L1 norm formulation contains non-differentiable absolute value terms. Due to these non-differentiable terms, we cannot directly apply Newton's Method to find stationary points from the optimization problem. 

\begin{figure}[h]
\centering
\subfloat[A]{\includegraphics[width=0.77in]{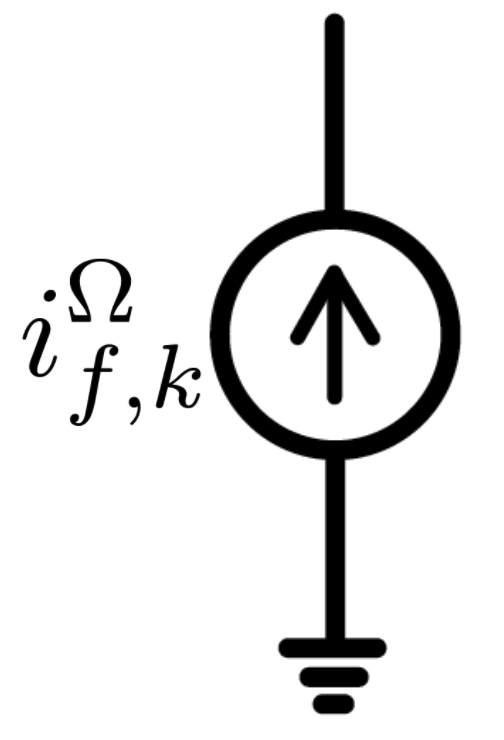}}
\label{sourceA}
\hfil
\subfloat[B]{\includegraphics[width=1.9in]{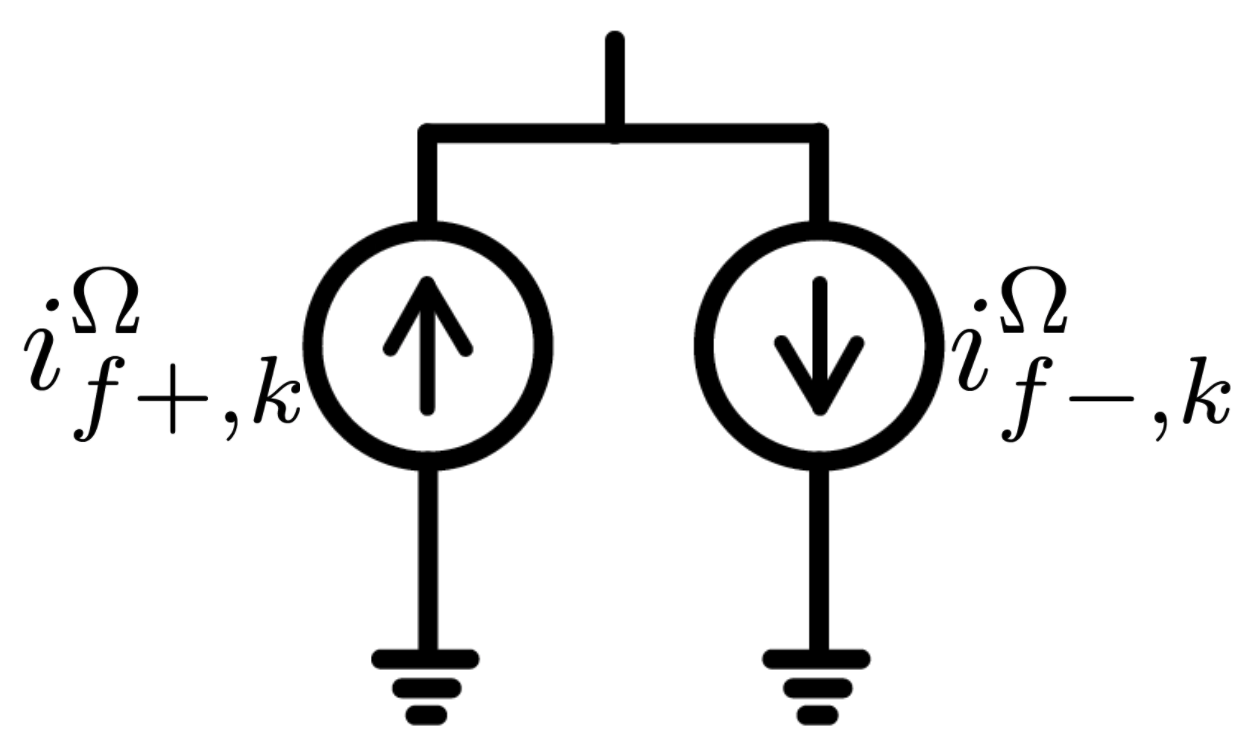}
\label{sourceB}}
\caption{[A] shows how infeasibility sources are modeled in the least squares formulation's circuit. [B] shows how infeasibility sources are modeled in the L1 norm formulation. In [B] the infeasibility source is represented as both an incoming and outgoing current source.}
\label{splitsource}
\end{figure}

To solve this problem, we can model infeasibility sources as two parallel current sources where one is flowing into the node and the other is flowing out of the node, as shown in Fig. \ref{splitsource}B, instead of the injection used in least squares shown in Fig. 
\squeezeupthree
\squeezeuptwo
\ref{splitsource}A. 
Mathematically we can represent these parallel sources by setting $i_{f,R}^{\Omega} = i_{f,R+}^{\Omega} - i_{f,R-}^{\Omega}$ and $i_{f,I}^{\Omega} = i_{f,I+}^{\Omega} - i_{f,I-}^{\Omega}$ and then constraining $i_{f,R+}^{\Omega}, i_{f,R-}^{\Omega}, i_{f,I+}^{\Omega}$ and $i_{f,I-}^{\Omega}$ to be positive values as is done in \cite{Singh}. Our new L1 objective function then becomes  $|i_{f,R+}^{\Omega}| + |i_{f,R-}^{\Omega}| + |i_{f,I+}^{\Omega}| + |i_{f,I-}^{\Omega}|$ which can be simplified to $i_{f,R+}^{\Omega} + i_{f,R-}^{\Omega} + i_{f,I+}^{\Omega} + i_{f,I-}^{\Omega}$ since we constrained each of these variables to be positive. Now each variable is differentiable and can be used within Newton's Method. This new formulation in shown in Equations \eqref{eqn: opt1_L1_improved} - \eqref{eqn: opt4_L1_improved}.

\begin{subequations}
\begin{align}
\min \quad & i_{f,R+}^{\Omega} + i_{f,R-}^{\Omega} + i_{f,I+}^{\Omega} + i_{f,I-}^{\Omega} \label{eqn: opt1_L1_improved} \\
\textrm{s.t.} \quad & GV^{\Omega}_R -  BV^{\Omega}_I +  I^{NL}_{R} (V^{\Omega}_R, V^{\Omega}_I, Q^{\Omega}) \nonumber \\ 
& \quad \quad \quad \quad \quad \quad \quad -(i_{f,R+}^{\Omega} - i_{f,R-}^{\Omega}) = 0 
\label{eqn: opt2_L1_improved}\\
& BV^{\Omega}_R + GV^{\Omega}_I + I^{NL}_{I}(V^{\Omega}_R, V^{\Omega}_I, Q^{\Omega}) \nonumber \\
&  \quad \quad \quad \quad \quad \quad \quad -(i_{f,I+}^{\Omega} - i_{f,I-}^{\Omega}) =  0 \label{eqn: opt3_L1_improved} \\
&  i_{f,R+}^{\Omega}, i_{f,R-}^{\Omega}, i_{f,I+}^{\Omega}, i_{f,I-}^{\Omega}  \geq  0 \label{eqn: opt4_L1_improved}
\end{align}
\end{subequations}
We then formulate the Lagrangian of the L1 norm formulation in \eqref{lagrangel1}. Since we now have inequality constraints, we add new variables $\mu_{R+}^\Omega, \mu_{R-}^\Omega, \mu_{I+}^\Omega,$ and $\mu_{I-}^\Omega$ as the inequality dual variable vectors. 

\begin{align}
     \mathcal{L}( i_{f,R+}^{\Omega},  & i_{f,R-}^{\Omega}, i_{f,I+}^{\Omega}, i_{f,I-}^{\Omega}, V_R^{\Omega}, V_I^{\Omega}, Q^{\Omega}, \lambda_R^{\Omega} , \lambda_I^{\Omega} ,  \nonumber \\
     \mu_{R+}^\Omega, & \mu_{R-}^\Omega, \mu_{I+}^\Omega, \mu_{I-}^\Omega ) = \nonumber \\
      &  i_{f,R+}^{\Omega} + i_{f,R-}^{\Omega} + i_{f,I+}^{\Omega} + i_{f,I-}^{\Omega} \nonumber \\
    + & (\lambda_R^{\Omega})^T \left(G V^{\Omega}_R - B V^{\Omega}_I + I^{NL}_{R}(V^{\Omega}_R, V^{\Omega}_I, Q^{\Omega})\right) \nonumber\\
    + & (\lambda_I^{\Omega})^T \left(B V^{\Omega}_R + G V^{\Omega}_I + I^{NL}_{I}(V^{\Omega}_R, V^{\Omega}_I, Q^\Omega)\right) \nonumber \\
    - & (\lambda_R^{\Omega})^T (i_{f,R+}^{\Omega} - i_{f,R-}^{\Omega}) - (\lambda_I^{\Omega})^T(i_{f,I+}^{\Omega} - i_{f,I-}^{\Omega}) \nonumber \\
    + & (\mu_{R+}^{\Omega})^T i_{f,R+}^{\Omega} + (\mu_{R-}^{\Omega})^Ti_{f,R-}^{\Omega} \nonumber \\
    + & (\mu_{I+}^{\Omega})^T i_{f,I+}^{\Omega} + (\mu_{I-}^{\Omega})^Ti_{f,I-}^{\Omega}
 \label{lagrangel1}
\end{align}

Having formulated the Lagrangian, we then derive the first-order optimality conditions for the KKT conditions which are shown in \eqref{L1statIR+} - \eqref{L1statII-}, \eqref{primfeas} - \eqref{compslack}. The stationarity constraints for $Q^{\Omega}, V^{\Omega}_R,$ and $V^{\Omega}_I$ remain the same as in  \eqref{L2statQ} - \eqref{L2statVI}. The new stationarity constraints for the $i_f^\Omega$ variables are shown in \eqref{L1statIR+} - \eqref{L1statII-}. 
\begin{subequations}
\begin{IEEEeqnarray}{lCl}
    \frac{\partial \mathcal{L}}{\partial i_{f,R+}^{\Omega}} & = & 1 - \lambda_R^{\Omega} + \mu_{R+}^{\Omega}= 0  \label{L1statIR+} \\
    \frac{\partial \mathcal{L}}{\partial i_{f,R-}^{\Omega}} & = & 1 + \lambda_R^{\Omega} + \mu_{R-}^{\Omega}= 0  \label{L1statIR-} \\
    \frac{\partial \mathcal{L}}{\partial i_{f,I+}^{\Omega}} & = & 1 - \lambda_I^{\Omega} + \mu_{I+}^{\Omega}= 0  \label{L1statII+} \\
    \frac{\partial \mathcal{L}}{\partial i_{f,I-}^{\Omega}} & = & 1 + \lambda_I^{\Omega} + \mu_{I-}^{\Omega}= 0  \label{L1statII-} 
\end{IEEEeqnarray}
\end{subequations}
The equality primal feasibility equations in \eqref{eqn: opt2} and \eqref{eqn: opt3} remain the same with just the appropriate substitution for $i^f_R$ and $i^f_I$. We now also have the inequality primal feasibility equations shown in \eqref{primfeas} and the dual feasibility equations shown in \eqref{dualfeas}. 
\begin{align}
    i_{f,R+}^{\Omega}, i_{f,R-}^{\Omega}, i_{f,I+}^{\Omega}, i_{f,I-}^{\Omega}  \geq  0
    \label{primfeas} \\
    \mu_{R+}^\Omega, \mu_{R-}^\Omega, \mu_{I+}^\Omega, \mu_{I-}^\Omega  \geq 0
    \label{dualfeas}
\end{align}
Because we introduced inequality constraints, we also add the complimentary slackness equations shown in  \eqref{compslack}. The new variable $\epsilon$ in \eqref{compslack} represents a small valued number that stems from the perturbation of KKT equations (see primal-dual interior point algorithm in Boyd \cite{Boyd}). 
\begin{align}
    (\mu_{R+}^{\Omega})^T i_{f,R+}^{\Omega} 
    & + (\mu_{R-}^{\Omega})^Ti_{f,R-}^{\Omega}
    \nonumber \\
    + & (\mu_{I+}^{\Omega})^T i_{f,I+}^{\Omega} + (\mu_{I-}^{\Omega})^Ti_{f,I-}^{\Omega} = \epsilon
    \label{compslack}
\end{align}

The sparsity achieved with the L1 norm formulation is discussed in Section V and is shown in Fig. \ref{fig_sim}B.

\section{Implementation of TPIA}
We developed this framework within our tool called SUGAR-distribution (SUGAR-D) \cite{code}, a version of which was used in the TPF framework in \cite{Amrit}. We built our own tool rather than using a commercial solver in order to best exploit the physics of the problem; ensure the real, rather than relaxed, problem is solved; and to incorporate circuit heuristics for improved convergence \cite{Amrit}. SUGAR-D uses ECF in conjunction with Newton's Method to solve TPF simulations. Algorithm 1 succinctly describes the TPIA methodology. 

The TPIA implementation includes but is not limited to the following electric power grid elements: capacitors; voltage regulators; transformers; switches; overhead, underground, and triplex transmission lines; inverter-based distributed generators (IBDGs); fuses; shunts; and loads. ECF element three-phase models draw on the physical nature of the component, so physical properties like different combinations of delta and wye connections, unbalance amongst phases, or the inclusion of the neutral wire can be incorporated without loss of generality.

With the L1 norm formulation we used a special limiting heuristic called diode limiting \cite{diode} to ensure convergence.  Diode limiting is used to mitigate large step sizes between Newton's Method iterations amongst inequality dual variables, ensuring that the dual variables remain within their feasible space. To accommodate the sensitivity of the L1 norm formulation to initial conditions, the formulation was also ``warm started" for certain testcases with the TPF solution to the testcase at a less stressed set of conditions. Since the solution to the L1 norm formulation will be the TPF solution if there are no infeasibility currents, the warm start enables the L1 norm formulation to have a reasonable initial starting point.

\begin{algorithm}
\begin{algorithmic}
	\caption{Infeasibility Analysis } 
	\label{alg:framework}
	\State {\bf procedure:}
	\State {\bf input} distribution network models, warm start (TPIA L1 norm only) or initial conditions, solver settings, Newton's Method tolerance ($\epsilon)$
	\State {\bf add} dual variables $(\mu, \lambda)$ and current sources $(i_f)$
	\State {\bf initialize} primal and dual variables $(\mu, \lambda)$, infeasibility currents $(i_f)$
	\State {\bf create} sparse matrix structure
	\State {\bf loop} through linear elements to create the linear admittance matrix $Y_{Lin}$
	\While{iteration count $k$ $<$ max iterations}
	\State{\bf loop} through nonlinear elements and update values
	\State{\bf solve} equations for unknown values at $k+1$ using Newton's Method
	\State{\bf limit} dual variables (TPIA L1 norm only only)
	\State \textbf{calculate} the error between iterations $k$ and $k+1$
	\If{$error > tolerance$}
        \State \textbf{update} variables dependent on  $k$ to $k+1$ values
    \Else
        \State \textbf{break} the loop and return converged states
    \EndIf
    \EndWhile
	\State {\bf Calculate} residual values from original nonlinear equations and determine if any KKT conditions are violated
\end{algorithmic}
\end{algorithm}

\begin{table*}[t]
 \centering
 \caption{Results from running varying testcases using GridLAB-D, SUGAR-D TPF, and the TPIA formulations.}
  \label{table_results}
  \begin{tabular}{|c|c|c|c|c|c|c|c|c|c|c|}
  \hline
    \multirow{2}{*}{} & \multirow{2}{*}{\textbf{GridLAB-D}} & \textbf{SUGAR-D} & \multicolumn{4}{|c}{\multirow{2}{*}{\textbf{TPIA Least Squares}}} & \multicolumn{4}{|c|}{\multirow{2}{*}{\textbf{TPIA L1 norm}}} \\ 
    & & \textbf{TPF} & \multicolumn{4}{|c|}{} & \multicolumn{4}{|c|}{} \\  \hline
    \multirow{2}{*}{Testcase} & \multirow{2}{*}{Converged} & \multirow{2}{*}{Converged} & Matrix & \multirow{2}{*}{Iterations} & \multirow{2}{*}{Time (s)} & \multirow{2}{*}{Nonzero $i_f$} & Matrix & \multirow{2}{*}{Iterations} & \multirow{2}{*}{Time (s)} & \multirow{2}{*}{Nonzero $i_f$} \\ 
    & & & Size & & &  & Size & & &  \\ \hline
    R4-12.47-1 & Yes & Yes & 42614  & 5 & 1.28 & 0 & 75254 & 5 & 2.08 & 0  \\ \hline
    R1-12.47-3\_OV & No & No & 2244  & 11 & 0.23 & 73 & 4032  & $92^1$ & $1.84^1$ &3 \\ \hline
    R2-25.00-1\_OV & No & No & 22612  & 20 & 2.5 & 766 & 39364  & $86^1$ & $15.43^1$ & 1  \\ \hline
    R3-12.47-3\_OV  & No & No & 142594  & 43 & 25.7 & 1504 & 251750 & $439^1$ & $443.4^1$ & 7 \\\hline
    R4-12-47-1\_EV  & No & No & 46294  & 28 & 7.033 &  324 & 83074 &  $20^1$& $6.14^1$& 1 \\\hline
    8500node\_OV & No & No & 116724 & 235 & 118.43 & 4760 & 219576 & $826^1$ & $711.25^1$ & 4  \\ \hline
  \end{tabular}
  \tablefootnote{}{These values were obtained with a warm start. The reported values represent the sum of the respective quantity from the L1 approach with the quantities it took to get the warm start values.}
 \end{table*}
\section{Experiments}
\subsection{Experimental Set-up}
To validate the efficacy of the TPIA approach, we run experiments on realistic but modified taxonomy feeders and a modified IEEE 8500node testcase. We simulate these networks using GridLAB-D, SUGAR-D TPF (described in \cite{Amrit}), and the least squares and L1 norm TPIA formulations. Table \ref{table_results} documents the cases that were used in these experiments. They are available on Github at https://github.com/emfoster/TPIA-Experiment-Testcases in the Gridlabd file format .glm. 

The taxonomy feeders R4-12.47-1, R1-12.47-3, R2-25.00-1, R3-12.47-3, and R4-12.47-1, respectively represent a small urban center, moderate suburban area, heavy suburban area, and a heavy urban area with a rural spur in four different climate regions \cite{PNNL}. R4-12.47-1 is used in its original form but the other networks are modified to make them infeasible by increasing their loading.  We identify these cases in Table \ref{table_results} with the addition of ``OV" to their titles. The testcase R4-12.47-1\_EV has electric vehicles added to it and is described in Section V.c. For the purpose of this paper, new infeasibility current sources are added at all nodes in each testcase with the exception of the slack bus.

All experiments were ran using an IDE on an i7 processor. Table \ref{table_results} documents the results of our experiments on these networks. The matrix size column represents the size of the inverted square matrix used in the Newton's Method after any rows and columns of zeros are removed. Time, iterations, and nonzero infeasibility currents are reported for both TPIA formulations. Our code was not optimized for speed and could be improved in future versions.

\color{black}
\squeezeupthree
\begin{figure*}[t]
\centering
\subfloat[A]{\includegraphics[width=3.38in]{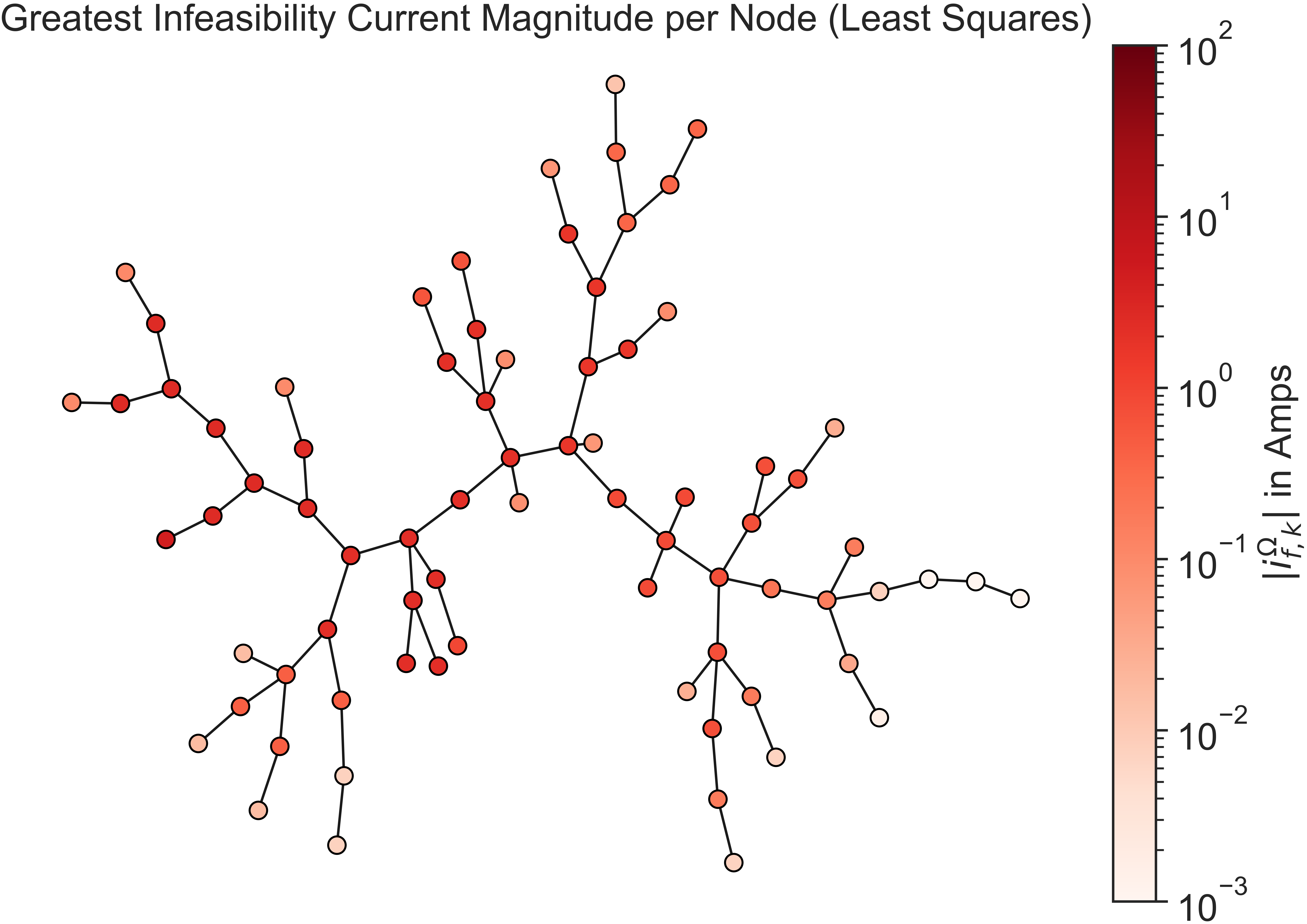}\label{fig_L2}}
\hfil
\subfloat[B]{\includegraphics[width=3.26in]{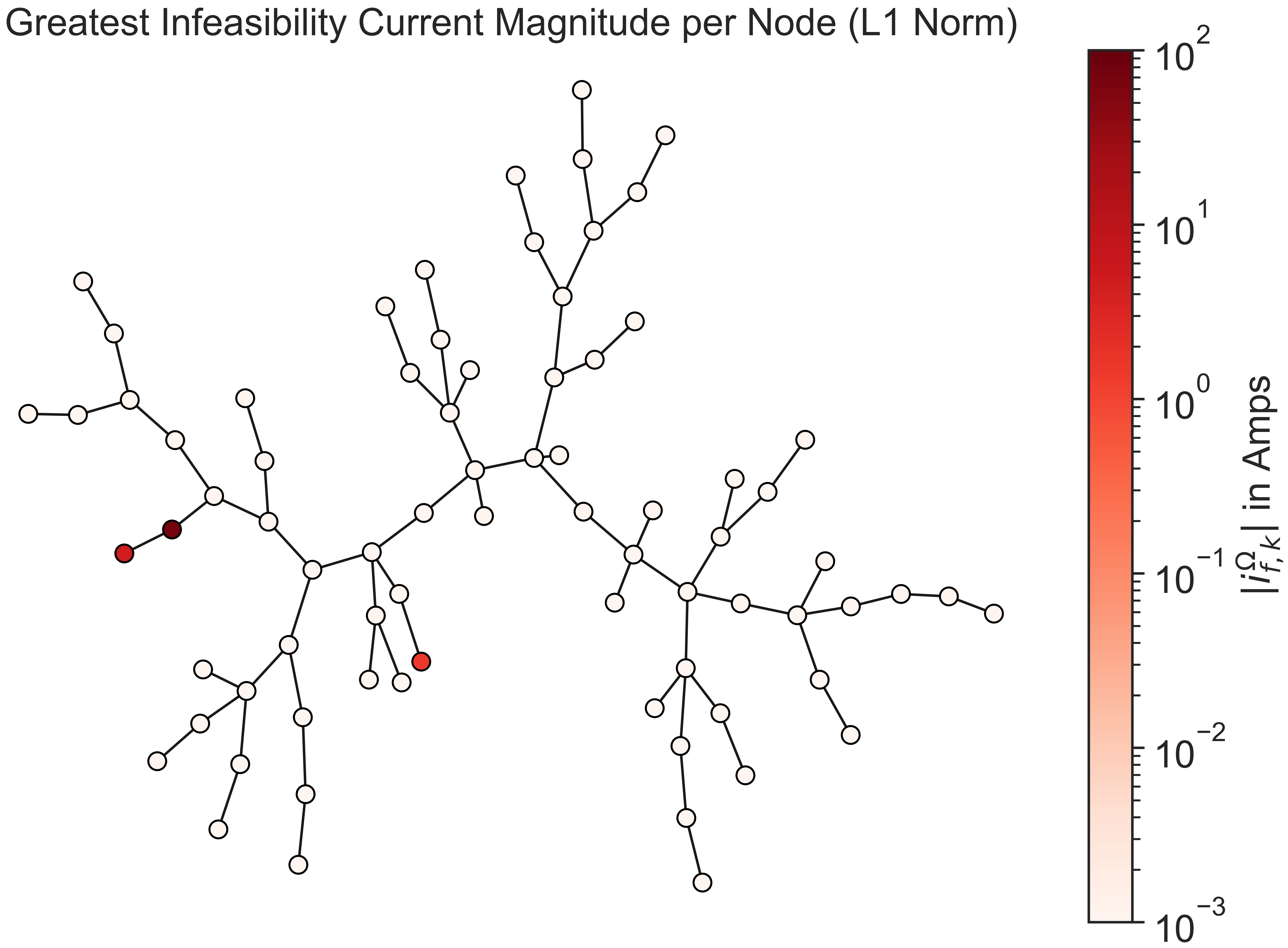}
\label{fig_L1}}
\caption{A graph representation of R1-12.47-3\_OV with a heatmap of the greatest infeasibility current magnitude amongst the phases at each node resulting from the least squares [A] and L1 norm [B] formulations. Currents are depicted on a logarithmic scale.}
\label{fig_sim}
\end{figure*}
\squeezeupthree


\subsection{Comparison of results between GridLAB-D, SUGAR-D TPF, and TPIA}
Table \ref{table_results} shows the results for each testcase obtained using GridLAB-D, SUGAR-D TPF, and the least squares and L1 norm TPIA formulations. GridLAB-D and SUGAR-D TPF only converge for the standard taxonomy feeder R4-12.47-1 which represents a traditional TPF case. On this testcase, we see that both TPIA approaches returned zero infeasibility currents which is what we would expect when TPF simulations converge. We also validated that the solution from the TPIA approaches matches with that from SUGAR-D TPF.

Neither GridLAB-D nor SUGAR-D TPF converged on any of the modified testcases that were designed to stress the system. Neither TPF simulations provided any insight into the physical reason for the failure. However, the two TPIA formulations converged to a solution with nonzero infeasibility currents for all modified testcases. The presence of non-zero infeasibility currents indicate that the problem in its original state is undispatchable.

For reporting purposes and to categorize arbitrarily small infeasibility current as zero, we create an infeasibility current threshold of $\epsilon_{If} = 1\text{e-}3 $ where values less than $\epsilon_{If}$ are set to zero. Table \ref{table_results} gives the number of maximal infeasibility current magnitudes above $\epsilon_{If}$. There is a stark difference in the number of nonzero infeasibility currents between the least squares and L1 norm formulations. In Fig. \ref{fig_sim}, the magnitude of the infeasibility currents at each phase was calculated and the maximum value amongst these phase infeasibility currents is shown at its respective node for R1-12.47-3\_OV. 

The least squares formulation is less sensitive to initial conditions than the L1 norm formulation and it did not require any additional heuristics to converge.  As shown in Table \ref{table_results}, least squares was faster and required fewer iterations than the L1 norm approach except for $R4-12.47-1\_EV$. Speed for both approaches could be improved through using additional circuits heuristics. For the L1 norm, a more optimal approach to finding a good warm start would also improve the time.

The infeasibility currents from TPIA least squares are spread throughout the network as shown in Fig. \ref{fig_sim}A, which makes it challenging for grid planners or operators to identify grid weaknesses. We can see from both Table \ref{table_results} and Fig.\ref{fig_sim}A that 73 of 76 nodes in the system have nonzero infeasibility currents for the least squares formulation. As noted in the Section III.b, infeasibility current sources can be added at only a subset of nodes, so for least squares, a better option may be for system planners with knowledge of the physical system to identify locations with potential for new or updated infrastructure and to conduct TPIA least squares over only that subset. In Fig. \ref{fig_sim}B, TPIA L1 norm had only 3 nodes with infeasibility currents. TPIA L1 norm is able to return sparse vectors of infeasibility currents, which in turn allows grid planners to isolate grid weakness and take realistic corrective action to make the network feasible. In a practical world setting, this implies that if the grid planner had access to the sparse L1 norm solution, they could add infrastructure to at most 3 locations in the network to make the case feasible; whereas, no such insight is available with the least squares solution since infeasibilities are spread across almost 96\% of the locations in Fig. \ref{fig_sim}.

\begin{figure}[t]
\centering
\includegraphics[scale = 0.35]{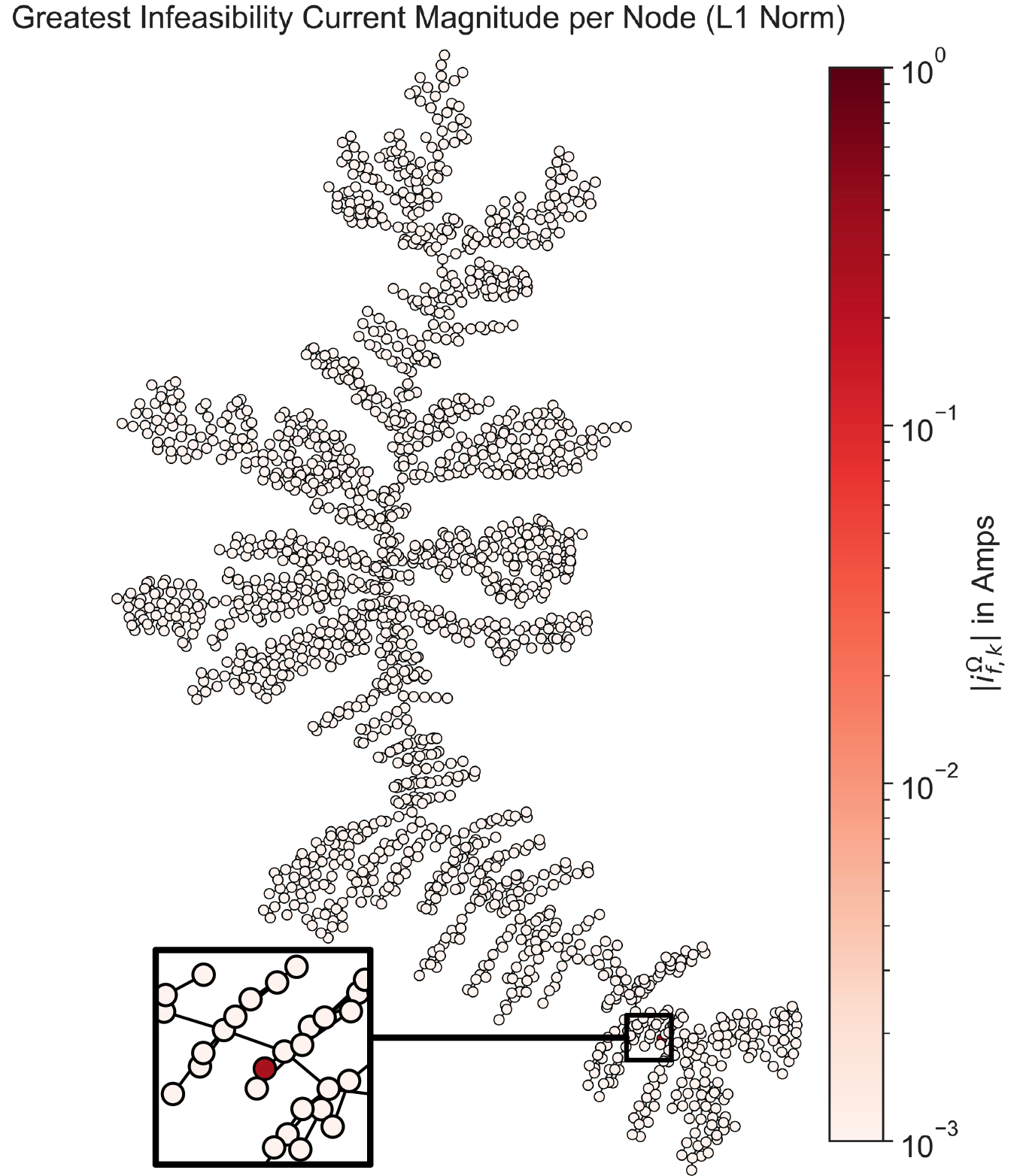}%
\label{fig_L1_BEV}
\caption{A graph representation of the testcase R4-12.47-1\_EV with a heatmap of the greatest infeasibility current magnitude amongst the phases at each node resulting from L1 norm formulation}
\label{fig_infeas}
\end{figure}
\squeezeupthree

\subsection{Use-case scenario on electric vehicle penetration and infrastructure additions}
To further illustrate how this research may be used in a planning scenario, we did a general study on the impact of adding battery electric vehicles (BEVs) to an overloaded version of testcase R4-12.47-1. From [17, Tab. DP04], 91.3\% of households in the U.S. have at least one vehicle available with 58.6\% of households having access to two or more vehicles. We conduct our study for the future scenario where internal combustion vehicles are phased out and we assume similar vehicle availability. We say that 90\% of households have at least one BEV available and that at some point 55\% of those vehicles are charging via a Level 2 charger conservatively at their maximum charge rate. We model the charging BEVs as simple constant PQ loads and add them at random triplex nodes in the system. We ran this testcase in GridLAB-D and SUGAR-D TPF and both failed to converge. We next used the L1 norm TPIA on the BEV network (R1-12.47-1\_EV) to localize areas of infeasibility in the system and found that there was one infeasible node in the network that was power deficient. The scenario is shown in Fig. \ref{fig_infeas} where we highlight the area of grid weakness. 

To demonstrate, how a grid planner could use this information to take corrective actions, we took the infeasibility currents and voltages for each phase for the node and then calculated the amount of missing power at each node. Then, to make the case feasible or dispatchable, we added a battery at the node rated at the calculated missing power. After adding the battery to the system, we were able to validate that the network is feasible by running the modified network using SUGAR-D TPF. We additionally ran both TPIA formulations and found that both return zero infeasibility currents with the added battery. The experiment also served as the validation of our choice of infeasibility current threshold $\epsilon_{If}$.  This is the value of infeasibility currents in the solution vector below which they are treated as inconsequential and as numerical noise.

The BEV use-case is a pilot example of how the TPIA methodology proposed within this paper can be used in practice. Other more sophisticated applications may significantly improve their performance with the use of this proposed methodology in the future.

\section{Conclusion}
 We developed a novel TPIA framework that converges for unbalanced networks that otherwise fail with traditional TPF. For both the least squares and L1 norm approaches, we show scalability of the method up to 8500 nodes. From the L1 norm formulation of TPIA we also extract critical information on where the cause of TPF failures might be for a given network. This proposed approach can be a critical tool for distribution planning studies moving forward. It can be used to study the transition of today's distribution grids to ones with a high-penetration of 
 DERs or inverter \-based technologies.  In future work, we plan to more rigorously evaluate the impact of new-age technologies on the distribution grid using this method.




%
\section*{Acknowledgment}
This work was sponsored in part by the  National Science Foundation under contract ECCS-1800812.

\end{document}